\documentclass[11pt]{amsart} 

\usepackage{amssymb,amsmath,latexsym,amsthm,amscd,eucal}
\usepackage{graphicx}
\usepackage{color}

\theoremstyle{plain}
\newtheorem{thm}{Theorem}[section]

\newtheorem{conj}[thm]{Conjecture}
\newtheorem{remark}[thm]{Remark}

\numberwithin{equation}{section}

\def\rank{\mathsf {rank}}

\def\ot{\otimes}

\newcommand{\half}{\frac{1}{2}}
\newcommand{\csf}{\mathsf{c}}
\newcommand{\g}{\mathfrak{g}}
\newcommand{\uf}{\mathsf{u}} 
\newcommand{\uni}{\mathfrak{n}}
\newcommand{\sms}{\mathsf{s}}

\newcommand{\As}{\mathsf{A}}
\newcommand{\N}{\mathsf{N}} 
 
\newcommand{\LL}{\mathsf{L}}

\newcommand{\M}{\mathsf{M}}

\newcommand{\Sf}{\mathsf{S}}

\newcommand{\UU}{\mathsf{U}}  
\newcommand{\V}{\mathsf{V}}
\newcommand{\W}{\mathsf{W}}  
\newcommand{\witt}{\mathsf{W}(1,\underline{1})}

\newcommand{\dimm}{\mathsf{dim\,}}

\begin{document}

\title[Problems in the Representation Theory of Modular Lie Algebras]
{Some Problems in the Representation Theory of Simple Modular Lie Algebras}

\author{Georgia Benkart}
\address{Department of Mathematics, University of Wisconsin-Madison, Madison, WI 53706-1325, USA} 
\email{benkart@math.wisc.edu} 

\author{J\"org Feldvoss}
\address{Department of Mathematics and Statistics, University of South Alabama, Mobile, AL 36688-0002, USA}
\email{jfeldvoss@southalabama.edu}

\thanks{\newline
\noindent 2010 \emph{Mathematics Subject Classification:} \ {Primary 17B50; Secondary 17B10, 17B20, 17B05}  \hfill  \smallskip \newline
\noindent \emph{Key words and phrases}: \  restricted Lie algebra, classical Lie algebra, Lie algebra of Cartan type, reduced enveloping
algebra, simple module, projective indecomposable module, Cartan matrix, block, representation type,
stable Auslander-Reiten quiver  \medskip  \newline  
$*$ The authors express their appreciation to Marina Avitabile, Thomas Weigel, and 
the Dipartimento di Matematica e Applicazioni at the Universit\`a degli Studi di
Milano-Bicocca for their warm hospitality during the \emph{Bicocca-Workshop on Lie Algebras 2013.}
They also thank Jim Humphreys for his comments on the first (arXiv) version of the paper.}

\vspace{-.8cm}
\dedicatory{To Helmut Strade with our best wishes on the occasion of his 70th birthday} 

\begin{abstract}
The finite-dimensional restricted simple Lie algebras of characteristic $p > 5$ are classical or
of Cartan type. The classical algebras are analogues of the simple complex Lie algebras and
have a well-advanced representation theory with important connections to Kazhdan-Lusztig
theory, quantum groups at roots of unity, and the representation theory of algebraic groups.
We survey progress that has been made towards developing a representation theory for the
restricted simple Cartan-type Lie algebras, discuss comparable results in the classical case,
formulate a couple of conjectures, and pose a dozen open problems for further study.
\end{abstract}

\maketitle

\date{}
 
\section{Introduction}

The finite-dimensional simple Lie algebras over an algebraically closed field of characteristic $p>3$
are either classical (analogues of the finite-dimensional simple complex Lie algebras), or belong to
one of four infinite families $\mathsf{W, S, H, K}$ of Cartan-type Lie algebras, or when $p=5$ are
Melikian algebras. A comprehensive exposition of this classification result can be found in the volumes
(\cite{St2}-\cite{St4}), see also the survey \cite{PrSt}. The Cartan-type Lie algebras are so named
because of their relation to the four families of infinite-dimensional Lie algebras that arose in Cartan's
pioneering study of pseudogroups in the early 1900's. The list of simple modular Lie algebras bears
a striking resemblance to the list of the simple complex Lie superalgebras, and many of the notions,
such as support varieties and complexity, that have been productive in studying representations
of simple Lie algebras of prime characteristic have been shown to yield important information about
representations of Lie superalgebras (cf.~\cite{Sr}, \cite{DuSr}, \cite{BoKN1}-\cite{BoKN4},
\cite{LeNZh}, \cite{BaKN}, \cite{Ba1}, \cite{Ba2}, and \cite{D}). Concepts that have proven useful
in the classification of the simple modular Lie algebras such as sandwich theory, for example, play
an essential role in Burnside-type problems in group theory. The study of the asymptotic behavior
of finite $p$-groups and their associated Engel Lie algebras has led in a very natural way to Lie
algebras of Cartan type. These coincidences hint at much deeper hidden connections linking these
topics. Understanding the representations of the simple Lie algebras of prime characteristic should
help to unravel some of these intriguing mysteries.

The representation theory of the classical Lie algebras has been the subject of many papers
that establish far-reaching connections with the representation theory of algebraic groups,
Kazhdan-Lusztig theory, and quantum groups at roots of unity (see for example,
\cite{FrP1}-\cite{FrP6}, \cite{Lu1}, \cite{Lu2}, \cite{H3}, \cite{J2}, \cite{J3}, and \cite{J6}).
There has been much progress in developing the representation theory of the simple Lie algebras
of Cartan type, however many questions remain. Our aim here is to discuss some of these
problems. Although we concentrate primarily on the restricted simple Lie algebras of classical
and Cartan types (as in \cite{BlWil}), many of the same questions could be raised for the 
restricted simple Melikian algebra, as well as for the minimal $p$-envelopes of the non-restricted
Cartan-type and Melikian Lie algebras. For some related problems,  we refer the reader to
\cite[Questions 1--3]{Sk4}.

\section{Background and notation}

Let  $\g$ be a finite-dimensional restricted Lie algebra over an algebraically closed field $\mathbf k$
of characteristic $p > 0$ with $p$-mapping $x \mapsto x^{[p]}$ for $x \in \g$. Any simple $\g$-module
$\M$ is finite dimensional and admits a character $\chi \in \g^\ast:=\mathsf{Hom}_{\mathbf  k}
(\g,\mathbf  k)$ such that the central element $x^p-x^{[p]}$ in the universal enveloping algebra
$\UU(\g)$ of $\g$  acts as the scalar $\chi(x)^p$ on $\M$  for all $x \in \g$. The simple $\g$-modules
with character $\chi$ are exactly the simple modules for the \emph{reduced  universal enveloping
algebra} $\uf(\g,\chi)$, which is $\UU(\g)$ factored by the ideal generated by the elements $x^p
-x^{[p]}-\chi(x)^p 1$ for all $x \in \g$. The dimension of $\uf(\g,\chi)$ is $p^{\dimm\g}$, and every
simple module for $\uf(\g,\chi)$ has dimension at most $p^{\lfloor\half\dimm\g\rfloor}$. If $\chi=0$,
then $\uf(\g):=\uf(\g,0)$ is the \emph{restricted universal enveloping algebra} of $\g$, and
$\uf(\g)$-modules $\M$ correspond to the restricted representations of $\g$ (where $x^{[p]}.m
=x^p.m$ for all $x \in \g$ and $m \in \M$). The extensions $\mathsf{Ext}^{\bullet}_{\UU(\g)}(\M,
\M')$ are trivial whenever $\M$ and $\M'$ are $\g$-modules having different characters. The dual
$\M^*$ of a $\uf(\g,\chi)$-module $\M$ is a module for $\uf(\g,-\chi)$, and the tensor product of
a $\uf(\g,\chi_1)$-module with a $\uf(\g,\chi_2)$-module is a module for $\uf(\g,\chi_1+\chi_2)$.
In particular, $\M^* \ot \M$ is a restricted $\g$-module.

Each character $\chi$ determines an alternating bilinear form on $\g$ defined by $\chi([x,y])$
for $x,y \in \g$, and the radical of the form,
\begin{equation}\label{eq:zchi}
\mathcal Z_{\g}(\chi) = \{z \in\g\mid\chi([z,y])=0\ \ \hbox{\rm for all} \ \ y \in \g\},
\end{equation}
often referred to as the \emph{centralizer of $\chi$}, has important connections with various invariants
of $\uf(\g,\chi)$-modules such as their support varieties.

If $\M$ is a finite-dimensional  indecomposable $\g$-module, there exists a character $\chi$ and a least
positive integer $r$ such that $(x^p-x^{[p]}-\chi(x)^p)^{p^{r-1}}.m = 0$ for all $x \in \g$ and $m\in\M$.
Let $\uf(\g,\chi)_r$  be the universal enveloping algebra $\UU(\g)$ factored by the ideal generated by
the elements $(x^p-x^{[p]}-\chi(x)^p)^{p^{r-1}}$ for all $x \in \g$. A finite-dimensional $\g$-module
$\M$ decomposes as a direct sum of indecomposable modules $\M = \bigoplus_j \M_j$, where each
$\M_j$ is a $\uf(\g, \chi_j)_{r_j}$-module for a uniquely determined character $\chi_j$ and a minimal
positive integer $r_j$ (see \cite[Thm.~5.2.6]{StF} and \cite[Thm.~2.2.2]{B}). The algebra $\uf(\g,\chi)_{r}$
is a cocommutative Hopf algebra when $\chi =0$, but  that need not be the case when $\chi$ is nonzero.
However,  $\uf(\g,\chi)_{r}$ is a Frobenius algebra for all $\chi$ and $r$ (see \cite[Cor.~5.4.3]{StF},
\cite[Prop.~1.2]{FrP4}, and \cite[Prop.~3.2.2]{B}). The comultiplication $\Delta: \uf(\g) \rightarrow
\uf(\g) \otimes \uf(\g)$, with $x\mapsto x\otimes 1+1 \otimes x$ for $x\in\g$, can be used to endow
$\uf(\g,\chi)$ with the structure of a left $\uf(\g)$-comodule algebra, and as a result, $\uf(\g,\chi)$ is a
right module algebra for the dual Hopf algebra $\uf(\g)^\ast$. A group-like element $\nu \in \uf(\g)^\ast$
defines a one-dimensional $\uf(\g)$-module $\mathbf k_\nu$ with $u.\alpha = \nu(u)\alpha$ for all $u \in
\uf(\g), \alpha \in \mathbf k$. The tensor product  $\M_\nu = \mathbf  k_\nu \otimes \M$ is then a $\uf(\g,
\chi)$-module for any $\uf(\g,\chi)$-module $\M$.  Thus, there is an action of the group of group-like elements
of $\uf(\g)^*$ on $\uf(\g,\chi)$-modules, which can be applied  in studying  the blocks of the algebra $\uf
(\g,\chi)$ (see \cite{F2}, \cite{Fe1}, and \cite{Fe3}).

An element $y$ of a restricted Lie algebra $\g$ is \emph{semisimple} (resp.~\emph{nilpotent}) if $y$
lies in the restricted subalgebra generated by $y^{[p]}$  (resp.~$y^{[p]^e} = 0$ for $e\gg 0$). Any
element $x\in\g$ has a Jordan-Chevalley decomposition $x=x_s+x_n$, where $x_s$ is semisimple, $x_n$
is  nilpotent, and $x_s$ and $x_n$ commute.  When $\g$ has a nondegenerate trace form, then a character
$\chi$ of $\g$ is said to be \emph{semisimple} (resp.~\emph{nilpotent}) if $\chi$ corresponds to a semisimple (resp.~nilpotent) element of $\g$ via the isomorphism from $\g^*$ onto $\g$ induced by the trace form.
In particular, when $\g$ is classical Lie algebra and $p$ is good for $\g$  (i.e., $p$ does not divide the
coefficient of any root relative to a basis of simple roots) and $\g$ is not isomorphic to $\mathfrak{sl}_{rp}
(\mathbf k)/\mathbf{k}\mathrm{I}$ for any $r\geq 1$, such a nondegenerate trace form exists (see
\cite{Ga}), so each character has a Jordan-Chevalley decomposition $\chi=\chi_s+\chi_n$ under these
assumptions. An alternate approach to decomposing $\chi$ can be found in \cite[Sec.~3]{KaW2}.

An element $t$ of a restricted Lie algebra $\g$  is \emph{toral} if $t^{[p]}=t$. A restricted subalgebra
$\mathfrak t$  of $\g$ is a \emph{toral subalgebra} (or a \emph{torus}) of $\g$ if the $p$-mapping is
invertible on $\mathfrak t$; equivalently, if $\mathfrak t$ is abelian and admits a basis consisting of toral
elements. The centralizer of a maximal torus (i.e., a torus not properly contained in another one) is a Cartan
subalgebra of $\g$, and every Cartan subalgebra is the centralizer of some maximal torus (see for example
\cite[Thm.~2.4.1]{StF}). For classical Lie algebras, maximal tori and Cartan subalgebras coincide, and all
have the same dimension. It was shown by Demuskin that if $\g$ is a restricted Lie algebra of Cartan
type, then all maximal tori of $\g$ have the same dimension and split into finitely many classes under the
action of the automorphism group of $\g$ (cf.~\cite[Sec.~7.5]{St2}). However, the Cartan subalgebras
of an arbitrary restricted Lie algebra may be non-conjugate under the automorphism group and may even
have different dimensions. For that reason, a natural notion for restricted Lie algebras is the \emph{maximal
toral rank},\
\begin{equation}\label{eq:mt}
\mathsf{mt}(\g):=\mathsf{max}\{\dimm\mathfrak{t}\mid\mathfrak{t}\ \hbox{\rm is a torus of\ $\g$}\}.
\end{equation}

Maximal tori are essential to understanding the structure and representation theory of restricted Lie
algebras. In particular, every finite-dimensional $\uf(\g,\chi)$-module $\M$ has a decomposition $\M
=\bigoplus_{\lambda\in\mathfrak t^*}\M_\lambda$ into weight spaces $\M_\lambda=\{m\in\M\mid
t.m=\lambda(t)m$ for all $t\in\mathfrak t\}$ relative to a maximal torus $\mathfrak t$.

The representation theory of the graded restricted Lie algebras is what is best understood at this
juncture. Each such algebra has a decomposition $\g=\bigoplus_{i=-q}^r\g_i$ into homogeneous
components with $[\g_i,\g_j]\subseteq\g_{i+j}$ if $-q\leq i+j\leq r$ and $[\g_i,\g_j]=0$ otherwise.
Let  $\g^{-}=\bigoplus_{i<0}\g_i$ and $\g^{+}=\bigoplus_{j > 0}\g_j$,  and assume $\chi$
is a character of $\g$ such that $\chi \vert_{\g^+} = 0$. Any simple $\g_0$-module $\LL_0$ with
character $\chi \vert_{\g_0}$, can be inflated to a module for the subalgebra $\mathfrak{h}:=\g_0
\oplus\g^+$ by having $\g^+$ act trivially. Then the induced module 
\begin{equation}
\mathsf{Z}(\LL_0,\chi):=\mathsf{Ind}_{\mathfrak h}^{\g}(\LL_0,\chi)=\uf(\g,\chi)
\ot_{\uf(\mathfrak{h},\chi\vert_{\mathfrak h})}\LL_0
\end{equation}
is a $\g$-module with character $\chi$.

According to \cite{S2}, when $\g$ is a graded restricted Lie algebra of Cartan type, any simple
$\uf(\g)$-module is isomorphic to a quotient of $\mathsf{Z}(\LL_0,0)$ for some simple $\g_0$-module
$\LL_0$. Moreover, the modules $\mathsf{Z}(\LL_0,\chi)$ for certain graded Cartan-type Lie algebras
and certain nonzero characters $\chi$ have been investigated in (\cite{Ch}, \cite{St1}, \cite{Ko1},
\cite{Ko2}, \cite{S1}-\cite{S3}, \cite{Ho1}-\cite{Ho5}, \cite{Hu1}, \cite{Hu2}, \cite{HoZ1}, \cite{HoZ2},
\cite{Z1}, \cite{Z2}, \cite{Sh2}, \cite{Sk1}, and \cite{J7}), where it was shown that in many instances
the induced modules $\mathsf{Z}(\LL_0,\chi)$ themselves are simple. For a graded Cartan-type Lie
algebra with the standard grading, the subalgebra $\g_0$ is a classical Lie algebra $\mathfrak{sl}_n(\mathbf k)$
or $\mathfrak{sp}_{2n}(\mathbf k)$ (or one of those algebras extended by a one-dimensional center). The simple
modules for $\uf(\g)$ have been described, but the description depends on knowing the simple modules
for $\uf(\g_0)$, which for large $p$ reduces to Lusztig's conjectures expressing the dimension of the
simple modules for classical Lie algebras in terms of Kazhdan-Lusztig polynomials. Understanding the
composition factors of the projective indecomposable modules of $\uf(\g)$ depends on knowing the
composition factors of the projective modules for $\uf(\g_0)$ (cf.~\cite[Thm.~3.1.5]{N} and
\cite[Thm.~4.1]{HoN2}).

For a graded restricted Lie algebra $\g=\bigoplus_{i=-q}^r\g_i$, the associated filtration is given by
$\g =\g_{(-q)}\supset\g_{(-q+1)}\supset\dots\supset\g_{(r)}\supsetneq 0$, where $\g_{(j)}:=
\bigoplus_{i=j}^r \g_i$ for all $-q \leq j\leq r$. The \emph{height} of a character $\chi$ of $\g$, defined by
\begin{equation}
\label{eq:hi}\mathsf{h}(\chi)=
\begin{cases}
\mathsf{min}\{-q\leq j\leq r\mid\chi\left(\g_{(j)}\right)=0\} & \hbox{\rm if}\,\,\,\chi(\g_r)=0,\\
r+1 & \hbox{\rm if}\,\,\,\chi(\g_r)\neq 0,
\end{cases}
\end{equation}
is a useful invariant, as the structure of a $\g$-module often is dictated by the height of its character.

A restricted Lie algebra $\g$ over an algebraically closed field $\mathbf k$ is said to have a \emph{triangular
decomposition} if $\g=\uni^-\oplus\mathfrak{t}\oplus\uni^+$, where $\mathfrak t$ is a maximal torus;
$\uni^-$ and $\uni^+$ are \emph{unipotent subalgebras} (i.e., restricted subalgebras such that
$(\uni^\pm)^{[p]^e}=0$ for $e\gg 0$); and $[\mathfrak t,\uni^\pm]\subseteq\uni^\pm$. Assume
$\chi \in\g^*$ is such that $\chi(\uni^+)=0$, and let $\mu\in\mathfrak t^*$ satisfy $\mu(t)^p-\mu(t)
=\chi(t)^p$ for all toral elements $t \in \mathfrak t$. Inflating the one-dimensional $\mathfrak t$-module
$\mathbf k_\mu$ determined by $\mu$ to the subalgebra $\mathfrak b^+:=\mathfrak t\oplus\uni^+$
and then inducing to $\uf(\g,\chi)$ yields the $\uf(\g,\chi)$-module $\mathsf{Z}(\mu,\chi):=
\mathsf{Ind}_{\mathfrak b^+}^\g(\mathbf k_\mu,\chi)$, the \emph{baby Verma module} determined
by $\mu$ and $\chi$. These modules feature prominently in the representation theory of the Lie algebras
with a triangular decomposition due to their universal property: every simple $\uf(\g,\chi)$-module is a
homomorphic image of $\mathsf{Z}(\mu,\chi)$ for some $\mu\in\mathfrak t^*$. The Lie algebra $\g$
of a connected reductive algebraic group $\mathsf{G}$ possesses such a triangular decomposition,
and there is always a character $\chi'$ in the $\mathsf{G}$-orbit of $\chi$ with the property that
$\chi'(\uni^+)=0$ (see \cite[1.4]{FrP4}, \cite[Sec.~10]{H3}, and \cite[B.2]{J6}), so there is no loss in
restricting considerations to such characters. Generalized baby Verma modules have been constructed
for the generalized reduced enveloping algebras $\uf(\g,\chi)_r$ of a Lie algebra $\g$ having a
triangular decomposition and have been used to study the simple and projective indecomposable
$\uf(\g,\chi)_r$-modules (see \cite{YShLi} and \cite{LiShY}).

Every finite-dimensional module of a finite-dimensional associative $\mathbf{k}$-algebra $\mathsf{A}$ has
a projective cover. The projective indecomposable $\mathsf{A}$-modules are in one-to-one correspondence
with the simple $\mathsf{A}$-modules $\mathsf{S}$ via the natural epimorphism $\mathsf{P(S)/Rad(P(S))}
\rightarrow \mathsf{S}$, which maps the projective cover $\mathsf{P(S)}$ modulo its radical $\mathsf{Rad(P(S))}$
onto $\mathsf S$. Let $\Xi$ index the simple $\mathsf{A}$-modules, and for an $\mathsf{A}$-module $\M$,
let $[\M]=\sum_{\mathsf{S}\in\Xi}\{\M:\mathsf{S}\}[\mathsf{S}]$ be the corresponding element in the
Grothendieck group, where $\{\M:\mathsf{S}\}$ is the nonnegative integer counting the multiplicity of
$\mathsf{S}$ as a composition factor of $\M$, and $[ \,\cdot\, ]$ denotes the isomorphism class. The
numbers $\{\mathsf{P(S):S'}\}$ for $\mathsf{S,S'\in\Xi}$ are the so-called \emph{Cartan invariants},
and the matrix with these numbers as its entries is the \emph{Cartan matrix}.

The finite-dimensional associative algebras $\As$ over an algebraically closed field can be separated
into three types (representation finite, tame, and wild), which  have the following description: (i) $\As$
is \emph{representation finite} or \emph{of finite representation type} if there are only finitely many
isomorphism classes of finite-dimensional indecomposable $\As$-modules. (ii) $\As$ is \emph{tame} or
\emph{of tame representation type} if it is not representation finite and if the isomorphism classes of
indecomposable $\As$-modules in any fixed dimension are almost all contained in a finite number of
one-parameter families. (iii) $\As$ is \emph{wild} or \emph{of wild representation type} if the category
of finite-dimensional $\As$-modules contains the category of finite-dimensional modules over the free
associative algebra in two generators. The classification of indecomposable objects (up to isomorphism)
for wild algebras is a well-known unsolvable problem, and so one is only able to classify the finite-dimensional indecomposable modules of algebras that are of finite or tame representation type.

\section{A dozen related problems}

Our primary focus will be on modules for restricted Lie algebras. All Lie algebras and their modules considered
here are finite dimensional over an algebraically closed field $\mathbf k$ of characteristic $p>0$.
\smallskip

\begin{itemize}
\item[]\textbf{Problem 1.} {\it For the reduced enveloping algebra $\uf(\g,\chi)$ of a restricted simple Lie
algebra $\g$ of Cartan type, determine the Cartan invariants $\{\mathsf{P(S)}:\mathsf{S'}\}$.}
\end{itemize}
\smallskip

For the Lie algebra $\g$ of a semisimple algebraic group, Humphreys \cite{H1} established a reciprocity
formula relating the Cartan matrix of its restricted enveloping algebra $\uf(\g)$ to the multiplicities
of the simple modules as composition factors of the baby Verma modules for $\uf(\g)$. This
formula inspired the well-known BGG reciprocity for semisimple complex Lie algebras as well as similar
reciprocity results  for Cartan-type Lie algebras by Nakano \cite{N} and Holmes and Nakano (\cite{HoN1},
\cite{HoN2}) and for generalized reduced enveloping algebras by Bendel \cite{B}. In \cite{J7}, Jantzen
generalized Humphreys' reciprocity law and showed for restricted simple Lie algebras of Cartan type that
the Cartan matrices corresponding to characters of height 0 are submatrices of the Cartan matrix for the
restricted enveloping algebra. 

The smallest of the restricted Cartan-type Lie algebras is the $p$-dimen\-sional Witt algebra $\witt$
with basis $\{e_n\mid n=-1,0,\dots, p-2\}$ and multiplication given by $[e_m, e_n]=(n-m)e_{m+n}$
(where $e_{m+n}=0$ if $m+n\not\in\{-1,0,\ldots,p-2\}$). The $p$-mapping on $\witt$ is specified by
$e_n^{[p]} = 0$ for $n \neq 0$ and $e_0^{[p]} = e_0$. Feldvoss and Nakano \cite{FeN} computed
the dimensions of the projective indecomposable modules $\mathsf{P}$ for all characters of $\witt$ as
well as the multiplicities $\{\mathsf{P}:\mathsf{S}\}$ of their composition factors.

When $\g$ is a restricted Lie algebra having a triangular decomposition $\g=\uni^-\oplus\mathfrak{t}
\oplus \uni^+$ such that the assumptions of \cite[Thm.~2.2]{HoN2} are satisfied, then $p^{\beta}$
divides all the Cartan invariants $\{\mathsf{P(S):S'}\}$ for $\uf(\g)$, where $\beta=\dimm\uni^+ -
\dimm\uni^- -\dimm\mathfrak{t}$. It follows for the Witt algebra $\witt$ that $p^{p-4}$ divides all
the Cartan invariants; while for the Hamiltonian Lie algebra $\mathsf{H}(2,\underline{1})$, $p^{p^2-10}$
divides all the Cartan invariants for its restricted enveloping algebra (see \cite{N}, \cite{HoN2},
\cite{Ho2}, \cite{ShJ}, and \cite{ShY}).
\smallskip

\begin{itemize}
\item[]\textbf{Problem 2.} \emph{If $\g(n)=\mathsf{X}(n,\underline{1})^{(2)}$, $n\geq 1$, are
the algebras in one of the families $\mathsf{X}=\mathsf{W,S,H}$, or $\mathsf{K}$ of restricted
simple Cartan-type Lie algebras, then using the fact that there is an embedding $\g(n)\hookrightarrow
\g(n')$ for $n<n'$, investigate the induction and restriction functors and the Grothendieck group
for $\varinjlim\g(n)$.}
\medskip

\item[]\textbf{Problem 3.} \emph{For a restricted simple Lie algebra of Cartan type $\g$, compute
the orbits of the characters under the action of the automorphism group and thereby determine
the isomorphism classes of the reduced enveloping algebras. Find a complete set of inequivalent
representatives for the orbits of the characters.}
\end{itemize}
\smallskip

For the Witt algebra $\witt$, a solution to Problem 3 has been achieved in \cite{FeN}. In contrast to
what occurs for classical Lie algebras (cf.~\cite[Sec.~2]{J5}), for $\witt$  there are infinitely many orbits of the
characters $\chi$ whose height $\mathsf{h}(\chi)$ is odd and satisfies $3\le\mathsf{h}(\chi)\le p-2$.
\smallskip

\begin{itemize}
\item[]\textbf{Problem 4.} {\it Investigate the $\mathsf{Ext}^1$ groups and $\mathsf{Ext}^1$
quivers for the restricted simple Cartan-type Lie algebras.}
\end{itemize}
\smallskip
 
For a restricted simple Lie algebra $\mathfrak{g}$ of classical type, Andersen \cite{A} proved that
the simple $\uf(\g)$-modules have no self-extensions (except for $\mathfrak{sp}_{2n}(\mathbf k)$
and $p=2$). However, in the extreme opposite case of a \emph{regular nilpotent character} $\chi$
(i.e., a nilpotent character $\chi$ for which $\dimm\mathcal{Z}_\g(\chi)=\mathsf{mt}(\g)$), there
are ``many" self-extensions of simple $\uf(\g,\chi)$-modules. In fact for such characters,
$\mathrm{Ext}_{\uf(\g,\chi)}^1(\Sf,\Sf)=0$ if and only if the simple $\uf(\g,\chi)$-module $\Sf$ is
projective (see \cite[Thm.~4.3\,(a)]{FrP4}).

Extensions between simple modules for the Witt algebra $\witt$ were determined in (\cite{BoNWi}
and \cite{R2}). For Cartan types $\W$ and $\mathsf{S}$, it was shown in \cite{LN} using a certain
infinite-dimensional Hopf algebra that simple modules have no self-extensions (see also \cite{CNPe}
for results on extensions of modules for $\witt$ and for the classical Lie algebras). For any
\emph{subregular nilpotent character} $\chi$ (i.e., a nilpotent character $\chi$ for which $\dimm
\mathcal{Z}_\g(\chi)=\mathsf{mt}(\g)+2$)  of a simple Lie algebra $\g$ of classical type, Jantzen
\cite{J4} determined the Loewy series of the baby Verma modules for $\uf(\g,\chi)$. Therefore,
under the assumptions in \cite{J4}, one obtains an explicit description of the simple modules and
the possible extensions between nonisomorphic simple modules.

It is important to compute the extensions $\mathsf{Ext}^1_{\uf(\g,\chi)}(\mathsf{S},
\mathsf{S}^\prime)$ of simple $\uf(\g,\chi)$-modules $\mathsf{S}$ and $\mathsf{S}^\prime$,
as this determines the \emph{Gabriel quiver} (= $\mathsf{Ext}^1$ quiver) of $\uf(\g,\chi)$
(without relations). If the cocycles are explicitly known, it would, at least in principle, be possible
to describe the structure of any finite-dimensional module by successively forming extensions
between its composition factors. Of course, this is too much to be expected in general, but even
partial answers can provide considerable information about the structure of representations.
Partial knowledge of the dimensions of the $\mathsf{Ext}^1$ spaces determines portions
of the Gabriel quiver and hence could provide necessary conditions for the representation
type of $\uf(\g,\chi)$. Moreover, the Gabriel quiver gives the following representation-theoretic
numerical  invariants:

\begin{enumerate}
\item[{$\bullet$}] the multiplicities of the second Loewy layer of the projective indecomposable
                             modules, and as a consequence, a lower bound for the Cartan invariants of
                             $\uf(\g,\chi)$;
\item[{$\bullet$}] the number of blocks of $\uf(\g,\chi)$ (which equals the number
                             of the connected components of the Gabriel quiver);
\item[{$\bullet$}] the number of isomorphism classes of simple modules in a given block of
                             $\uf(\g,\chi)$.
\end{enumerate}

Assume that $\g$ is a classical Lie algebra and $\chi$ is a nilpotent character of $\g$, so that
$\chi$ vanishes on some Borel subalgebra $\mathfrak{b}$ of $\g$. The Weyl group $\W$ of
$\g$ acts on the set of weights $\Lambda$ with respect to a maximal torus $\mathfrak{t}$
of $\mathfrak b$ via the dot action $w\cdot\lambda:=w(\lambda+\rho)-\rho$ (where $\rho$
denotes the sum of the fundamental weights). Brown and Gordon \cite[Thm.~3.18]{BrG} have
shown that the blocks of $\uf(\g,\chi)$ are in one-to-one correspondence with the orbits of
$\W$ on $\Lambda$, thereby confirming a conjecture of Humphreys (see \cite[Sec.~18]{H3}).
As a consequence, all reduced enveloping algebras $\uf(\g,\chi)$ for a nilpotent character $\chi$
of a classical Lie algebra $\g$ have the same number of blocks.

Humphreys has also speculated that for a nilpotent character $\chi$, the number of isomorphism
classes of simple modules in a given block of $\uf(\g,\chi)$ is at most the order of the Weyl group
of $\g$. This is true when $\chi$  is of \emph{standard Levi form}, that is to say, the restriction of
$\chi$ to some Levi subalgebra of $\g$ is a regular nilpotent character.
When $p$ is bigger than the Coxeter number, Corollary 5.4.3 and 
Theorem 7.1.1 of  \cite{BeMRu} give a precise count of the number of isomorphism 
classes of simple modules in a block of $\uf(\g,\chi)$ for any nilpotent character $\chi$ in 
terms of the Euler characteristics of the Springer fibers.

Holmes and Nakano \cite{HoN2} developed criteria for when a restricted Lie algebra has a unique
block and used them to prove that the restricted enveloping algebra $\uf(\g)$ of a restricted simple
Lie algebra $\g$ of Cartan type has a unique block.
\smallskip

\begin{itemize}
\item[]\textbf{Problem 5.} {\it Investigate the structure of the reduced enveloping algebras
$\uf(\g,\chi)$ for $\g$ a restricted simple Lie algebra of Cartan type {\rm (}which is most complicated
when $\chi=0${\rm )}. Give necessary and suffi\-cient conditions for $\uf(\g,\chi)$ to be semisimple.
Determine the Jacobson radical of $\uf(\g,\chi)$ {\rm (}which most likely is largest when $\chi=0${\rm)}.}
\end{itemize}
\smallskip

The set $\mathcal S_\g=\{\chi\in\g^*\mid\uf(\g,\chi)$ is semisimple$\}$ is Zariski open in $\g^*$
(see \cite[Thm.~4.5\,(3)]{F1} or \cite[Prop.~4.2\,(2)]{PrSk}). As a result, when $\mathcal S_\g$
is nonempty for a restricted Lie algebra $\g$, then $\mathcal S_\g$ is dense in $\g^*$, and thus
$\uf(\g,\chi)$ is semisimple for ``most'' characters $\chi$.

According to \cite[Cor.~3.6]{FrP4} (see also \cite[Sec.~14]{H3}), when $\g$ is classical and $p$
is good for $\g$, $\uf(\g,\chi)$ is semisimple if and only if the character $\chi$ is regular semisimple
(that is, $\chi$ is semisimple and $\dimm\mathcal{Z}_\g(\chi)=\mathsf{mt}(\g)$, or equivalently,
$\chi(t)\ne 0$ for every toral basis element $t$ of a maximal torus of $\g$). Because a classical
Lie algebra $\g$ always has a regular semisimple character, $\mathcal S_\g$ is nonempty, and
therefore $\uf(\g,\chi)$ is semisimple for ``most'' characters $\chi$ of a classical Lie algebra.

Premet and Skryabin \cite{PrSk} have conjectured that an arbitrary restricted Lie algebra $\g$ is
\emph{generically semisimple} (i.e., $\mathcal S_\g\neq\emptyset$) if and only if there exists
$\xi\in\g^*$ such that the centralizer $\mathcal Z_\g(\xi)$ (cf.~\eqref{eq:zchi}) is a toral subalgebra
of $\g$. They have established one of the directions: the existence of such a $\xi$ implies that
$\g$ is generically semisimple. The converse has recently been shown for solvable restricted Lie
algebras when $p>2$ (see \cite[Thm.~2.3]{Sk3}). Furthermore, Skryabin has conjectured that
$\uf(\g,\chi)$ is a simple algebra if and only if $\mathcal Z_\g(\chi)=0$. In \cite[Thm.~1\,(i)]{Sk2},
this conjecture is confirmed for solvable restricted Lie algebras when $p>2$, and in \cite[Thm.~2]{Sk2},
Skryabin has verified that it also holds under the assumptions that $\g$ is a Frobenius Lie algebra
(i.e., $\mathcal Z_\g(\xi)=0$ for some $\xi\in\g^*$) and all the inner derivations of $\g$ belong to
the Lie algebra of the automorphism group of $\g$.
\smallskip

\begin{itemize}
\item[]\textbf{Problem 6.} {\it For $\g$ a restricted simple Lie algebra, determine when $\uf(\g,\chi)$
is Morita equivalent to $\uf(\g',\chi')$ for some semisimple Lie algebra $\g'$ of smaller dimension {\rm (}
\hspace{-.2cm}or possibly to $\uf(\g',\chi')$ tensored with a commutative associative semisimple
algebra{\rm )}}.
\end{itemize}
\smallskip

 For a classical Lie algebra $\g$ and $p$ good for $\g$,   such results are possible (see for example,
 \cite[Thm.~2]{KaW1}, \cite[Thm.~3.2 and Thm.~8.5]{FrP4}, or \cite[Thm.~B.8]{J6}), 
and they can be used to reduce
the representation theory of $\uf(\g,\chi)$ for an arbitrary character $\chi$ to the representation
theory of $\uf(\g',\chi\vert_{\g'})$ for some reductive restricted subalgebra $\g'$ of $\g$ such that
$\chi\vert_{\g'}$ is nilpotent  as in  (\cite[Sec.~15]{H3} and \cite[B.8 and B.9]{J6}).

A Lie algebra $\g$ can be regarded as a $\g$-module under the adjoint action $\mathsf{ad}_x(y)
= [x,y]$. Hence, the symmetric algebra $\mathsf{S}(\g)$ and the quotient algebra $\sms(\g,\chi)$,
which is $\mathsf{S}(\g)$ modulo the $\g$-invariant  ideal generated by the elements $(x-\chi(x)1)^p$
for all $x\in\g$, are also $\g$-modules. Moreover, $\dimm\sms(\g,\chi)=\dimm\uf(\g,\chi)=p^{\dimm\g}$.
The Lie algebra $\g$ is said to have the {\it Richardson property} if there exists a faithful simple
$\g$-module $\mathsf{Y}$ with associated representation $\varrho:\g\rightarrow\mathfrak{gl}
(\mathsf{Y})$ such that $\mathfrak{gl}(\mathsf{Y})=\varrho(\g)\oplus\mathsf{R}$, where
$\mathsf{R}$ is a subspace of $\mathfrak{gl}(\mathsf{Y})$ with $[\varrho(\g),\mathsf{R}]\subseteq
\mathsf{R}$. For a Lie algebra $\g$ with the Richardson property and $\chi\in\g^*$, it follows
from \cite[Prop.~2.3]{Pr3} that ${\sms}(\g,\chi)$ and $\uf(\g,\chi)$ are isomorphic $\g$-modules.
Lemma 3.7 of \cite{Pr3} shows that the Witt algebra $\witt$ does not have the Richardson property.
In \cite{Pr3}, Premet speculates that the Richardson property characterizes the class of restricted
Lie algebras  $\g$ of the form $\g=\mathsf{Lie (G)}$, where $\mathsf{G}$ is a reductive
$\mathbf{k}$-group such that $p$ is good for $\g$ and $\g$ is not isomorphic to $\mathfrak{sl}_{rp}
(\mathbf k)/\mathbf{k}\mathrm{I}$ for any $r\geq 1$.
\smallskip

\begin{itemize}
\item[]\textbf{Problem 7.} \it{Assume that $\g$ is an arbitrary restricted Lie algebra.}
\begin{itemize}
\item[{\rm (a)}] \it{For $\V$ a simple $\g$-module, investigate  the structure of the tensor $\g$-module
$\mathsf{T}(\V)=\bigoplus_{k \geq 0}\V^{\otimes k}$,\, the symmetric tensors $\mathsf{S}(\V)=
\bigoplus_{k\geq 0}\mathsf{S}^k(\V)$, and the exterior tensors $\Lambda(\V)=\bigoplus_{k\geq 0}
\Lambda^k(\V)$. Determine the $\g$-module invariants {\rm (}i.e., the sum of the trivial $\g$-modules{\rm )}
in these modules.}
\item[{\rm (b)}] Determine the structure of the $\g$-module $\sms(\g,\chi)$ for $\chi\in\g^\ast$.
\item[{\rm (c)}] Prove/disprove Premet's conjecture on the Richardson property.
\end{itemize}
\end{itemize}
\smallskip

For a finite-dimensional Hopf algebra $\mathsf{A}$ and an $\mathsf{A}$-module $\V$, if the annihilator $\mathsf{Ann}_\mathsf{A}(\V)$ contains no nonzero Hopf ideal, then every projective indecomposable
$\mathsf{A}$-module is isomorphic to a direct summand of $\V^{\ot n}$ for some $n$ such that $1\leq
n\leq\mathsf{d}(\mathsf{A})$, where $\mathsf{d}(\mathsf{A})$ is the maximal length of a strictly
descending chain of proper two-sided ideals in $\mathsf{A}$ (see \cite{FeKl}). When these hypotheses
are satisfied, the tensor powers of a well-chosen module $\V$ can be used to construct all the projective
indecomposable $\mathsf{A}$-modules. When $\g$ is a restricted simple Lie algebra, every projective
indecomposable $\uf(\g)$-module is a direct summand of a suitable tensor power of the adjoint module
(see \cite[Cor.~1]{Fe2}).

Assume $\g$ is a restricted Lie algebra over $\mathbf k$, and $\M$ is a $\g$-module with character
$\chi$. Let $\vert(\g,\chi)\vert_\M$ be the affine homogeneous variety associated with the annihilator
of the Noetherian commutative algebra $\mathsf{H}^{2\bullet}(\g,\mathbf k)$ on $\mathsf{H}^{\bullet}
(\uf(\g),\M^*\ot\M)$ (see \cite{FrP1}). There exists a natural finite morphism of varieties, $\vert(\g,0)
\vert_\mathbf k\rightarrow\g^{(-1)}$, where $\g^{(-1)}$  is the vector space $\g$ with the scalars
$\alpha \in \mathbf k$ acting by multiplication by $\alpha^p$. By (\cite{J1}, \cite{FrP3}, and
\cite[Thm.~6.4]{FrP4}) the image of the variety $\vert(\g,\chi)\vert_\M$ under this morphism is the
{\it support variety} $\mathcal V_\g(\M)$, which can be characterized by
\begin{equation}
\mathcal V_\g(\M)\,=\,\{x\in\g\mid x^{[p]}=0\ \ \hbox{\rm and \ $\M\vert_{\langle x\rangle_p}$\ \ is not free}\}
\cup\{0\},
\end{equation}
where $\langle x\rangle_p$ is the restricted subalgebra of $\g$ generated by $x$. The theory of
support varieties for restricted Lie algebras has been developed by Friedlander and Parshall in a series
\cite{FrP2}-\cite{FrP6} of foundational papers (see also \cite{J1} and the survey \cite{J2} for results
on the cohomology and support varieties of restricted Lie algebras). In \cite[Prop.~6.2]{FrP4}, it is
shown that $\mathcal V_\g(\M)=0$ if and only if $\M$ is a projective (or injective) module for $\uf(\g,\chi)$.
Moreover, the following holds: 

\begin{thm}{\rm (See \cite[Prop.~2.1a)]{FrP3}, \cite[Prop.~7.1(a) and Cor.~7.2]{FrP4}, and
\cite[3.4(3)]{J2}.)} Assume $\g$ is a restricted Lie algebra, $\mathfrak h$ is a restricted subalgebra
of $\g$, and $\chi\in\g^*$. Let $\M$ be a $\uf(\g,\chi)$-module. Then
\begin{itemize}
\item[{\rm (a)}] $\mathcal V_\mathfrak h(\M\vert_\mathfrak h)=\mathcal V_\g(\M)\cap\mathfrak h$.
\item[{\rm (b)}] Any closed, conical subvariety of $\mathcal V_\g(\M)$ is of the form $\mathcal V_\g
(\mathsf N)$ for some $\g$-module $\mathsf{N}$ with character $\chi$.
\end{itemize}
\end{thm} 

Assume now that $\g$ is the Lie algebra of $\mathsf{G}$, a connected and simply connected reductive
group, and $\M$ is a $\uf(\g,\chi)$-module. Then $\mathcal V_\g(\M)\subseteq\mathcal Z_\g(\chi)$ by
\cite[Corrigenda and Addenda, Thm.~1.1${}'$]{Pr2}. Premet remarks that it is plausible this result holds
for a wider class of restricted Lie algebras; however, an explicit example is provided in \cite{Pr2} to
demonstrate that the result  is not true in general. Premet conjectured in \cite{Pr2} and showed in a
sequel paper \cite{Pr3}, that when $p$ is good for $\g=\mathsf{Lie(G)}$, then for every $\chi \in \g^*$ there exists
a $\uf(\g,\chi)$-module $\M$ such that $\mathcal V_\g(\M)=\mathcal{N}_p(\g)\cap\mathcal{Z}_\g
(\chi)$, where $\mathcal{N}_p(\g)=\{x\in\g\mid x^{[p]}=0\}$. This enabled him to determine for
which $\chi \in \g^*$ the algebra $\uf(\g,\chi)$ has finitely many indecomposable modules up to
isomorphism and to obtain a result towards determining which algebras $\uf(\g,\chi)$ are of tame type.
Subsequent papers (\cite{NPo}, \cite{BrG}, and \cite{GPr}) have expanded upon this theme.

We have the following characterizations of when $\uf(\g)$ has finite representation type.

\begin{thm}{\rm (Cf. \cite[Thm.~3]{PfV}, \cite[Thm.~2.4]{FeSt}, \cite[Thm.~4.3]{F1},
and \cite[Thm.~2.7]{FV1}.)}\label{T:finite}
The following four statements are equivalent for a restricted Lie algebra $\g$:
\begin{itemize}
\item $\uf(\g)$  has finite representation type.
\item There exists a toral element $t\in \g$  and a nilpotent element $x \in \g$ such that $\g=
\mathbf{k}t \oplus \mathsf{N}(\g)$,  where $\mathsf{N}(\g)=\mathsf{T}(\g)\oplus\langle x
\rangle_p$; $\mathsf{N}(\g)$ is the largest nilpotent ideal of $\g$; $\mathsf{T}(\g)$ is the
largest restricted ideal of $\g$ that is a torus; and $\langle x\rangle_p$ is the restricted
subalgebra generated by $x$.
\item $\dimm\mathcal{V}_\g(\mathbf k)\leq 1$.
\item $\uf(\g,\chi)$ has finite representation type for every $\chi\in\g^*$. 
\end{itemize}
\end{thm} 

The set $\mathcal F_\g=\{\chi\in\g^*\mid\uf(\g,\chi)$ has finite representation type$\}$ is
Zariski open in $\g^*$ (see \cite[Thm.~4.5\,(2)]{F1}). So when $\mathcal F_\g$ is nonempty
for a restricted Lie algebra $\g$, then $\mathcal F_\g$ is dense in $\g^*$, and thus $\uf(\g,
\chi)$ has finite representation type for ``most'' characters $\chi$.

A generalization of Theorem \ref{T:finite} to the restricted enveloping algebras of tame type has
been obtained for $p>2$ in (\cite{V}, \cite{FV2}, \cite{FV3}, \cite{FSko}, and \cite{F6}).
\smallskip

To understand the module category of an associative algebra $\As$, it is essential to know not
only the indecomposable objects  but  also the \emph{irreducible morphisms}  between them and
the corresponding \emph{Auslander-Reiten sequences}. A good starting point is to obtain information
about the connected components of the \emph{stable Auslander-Reiten quiver}. The vertices of
the stable Auslander-Reiten quiver of $\As$ are the isomorphism classes of non-projective
indecomposable $\As$-modules, and the arrows are given by the irreducible morphisms. According
to Riedtmann's structure theorem, these connected components  can  be described by certain directed
trees (and admissible subgroups of their automorphism groups). In the case of modular group algebras,
the underlying graphs of the connected components are finite, infinite, or affine Dynkin diagrams. 
Erdmann \cite{E} showed that the same is true for the restricted enveloping algebra of every
restricted Lie algebra $\g$, and in fact her proof is valid for any reduced enveloping algebra $\uf
(\g,\chi)$ (see \cite{F1}). It remains to determine which Dynkin diagrams actually occur. In the group
case, this is known due to the work of Okuyama, Erdmann, Bessenrodt, and others. However, many
of the methods used for group algebras are not available for Lie algebras.
\smallskip  

\begin{itemize}
\item[]\textbf{Problem 8.} \emph{Let $\g$ be a restricted simple Lie algebra of Cartan type,
and let $\chi\in\g^*$.}
\begin{itemize}
\item[{\rm (a)}] \emph{Determine the representation type {\rm (}finite, tame, wild{\rm )} of the blocks of $\uf(\g,\chi)$.}
\item[{\rm (b)}] \emph{Determine the connected components of the stable Auslander-Reiten quiver of $\uf(\g,\chi)$.}
\item[{\rm (c)}] \emph{Identify the position of the simple $\uf(\g,\chi)$-modules in the
                          stable Aus\-lander-Reiten quiver of $\uf(\g,\chi)$}.
\end{itemize}
\end{itemize}
\smallskip  

Part (a) is known for the Witt algebra $\g=\witt$ when $p\ge 7$.  More specifically, if
$\mathsf{h}(\chi)\leq p-4$, then $\uf(\g,\chi)$ is wild. If $\mathsf{h}(\chi)\ge p-3$,
then $\uf(\g,\chi)$  is a Nakayama algebra, that is to say, the projective indecomposable
modules have a unique composition series (see \cite[Thm.~5.2]{F2} for details). However,
when $p=5$, $\uf(\g,\chi)$ is wild if $\mathsf{h}(\chi)=-1,0$, and $\uf(\g,\chi)$ is a Nakayama
algebra if $\mathsf{h}(\chi)\ge 2$ (see \cite[Sec.~4]{FeN}), but the representation type
of $\uf(\g,\chi)$ is not known for characters $\chi$ of height 1 (see \cite[Sec.~4]{FeN} and
\cite[Thm.~5.3.1]{R1}).

By \cite[Thm.~6.6.2]{B}, for an arbitrary  restricted Lie algebra $\g$, the representation type
of the generalized reduced enveloping algebra $\uf(\g,\chi)_{r}$ for $r>1$ is wild when $\dimm
\g\geq 3$, and if $\uf(\g,\chi)_{r}$ is of tame representation type, then $\dimm\g=2$.

Results on the stable Auslander-Reiten quiver can be found for the Witt algebra $\witt$ in
\cite[Ex.~5.3]{F3}; for more general Cartan-type Lie algebras in \cite[Cor.~3.5]{F4}; and for
Lie algebras of reductive groups in \cite[Thm.~5.2 and Thm.~5.5]{F5}. When two vertices $[\M]$
and $[\N]$ belong to the same connected component $\Theta$ of the stable Auslander-Reiten
quiver, they have equal support varieties (see \cite[Lem.~5.2]{F1}), and one can then set
$\mathcal V_\g(\Theta):=\mathcal V_\g(\M)$ for any vertex $[\M]$ of $\Theta$. When $\dimm
\mathcal{V}_\g(\Theta) \geq 3$ for  a restricted Lie algebra $\g$, the component $\Theta$ has
tree class $\mathsf{A}_\infty$, but the converse of that statement fails to hold (see \cite{F4}
for details). Rickard \cite{Ri} (see also \cite{F7} for a correction) has shown that if $\dimm
\mathcal V_\g(\M)\geq 3$ for a $\uf(\g,\chi)$-module $\M$, then $\uf(\g,\chi)$ is of wild representation
type. That result is used in the proof of \cite[Cor.~3.4]{F4} to deduce that if $\g$ is a classical Lie
algebra not isomorphic to $\mathfrak{sl}_2(\mathbf k)$, then any block of $\uf(\g)$ not associated to the
Steinberg module is necessarily wild. (It is well known that any block of $\uf(\mathfrak{sl}_2(\mathbf k))$
not associated to the Steinberg module is tame, see \cite{Fi}.)

Farnsteiner and R\"ohrle \cite[Thm.~7.1]{FR1} have shown that when $p\ge
5$, the non-projective baby Verma modules for the restricted enveloping
algebra $\uf(\g)$ of a classical Lie algebra $\g$ belong to a stable
Auslander-Reiten component of tree class $\mathsf{A}_\infty$. They either are
periodic or have exactly one predecessor in the stable Auslander-Reiten quiver.
By \cite{FR2}, the non-projective simple $\uf(\g)$-modules are contained in a stable Auslander-Reiten
component of tree class $\tilde{\mathsf{A}}_{12}$ or $\mathsf{A}_\infty$. In the
$\mathsf{A}_\infty$-case, the simple module $\Sf$ has exactly one predecessor, and the middle
term of the Auslander-Reiten sequence terminating in $\Sf$ is indecomposable. As a consequence,
the \emph{heart} of the projective cover $\mathsf{P(S)}$ of $\Sf$ (i.e., the radical of $\mathsf{P(S)}$
modulo its socle) is either indecomposable or a direct sum of two copies of a simple $\uf(\g)$-module
not isomorphic to $\Sf$ (see \cite[Cor.~2.4]{FR2}).
\medskip  \smallskip

\begin{itemize}
\item[]\textbf{Problem 9.} \emph{Assume that $\g$ is an arbitrary restricted Lie algebra and $\chi
\in\g^*$.}
\begin{itemize}
\item[{\rm (a)}] \emph{Let $\mathsf{N}$ be a $\uf(\mathfrak h,\chi\vert_{\mathfrak h})$-module
for a  restricted subalgebra $\mathfrak h$ of $\g$, and consider the induced module $\M=
\mathsf{Ind}_\mathfrak h^\g(\mathsf{N,\chi})$. Then $\mathcal V_\g(\M)=\mathcal V_{\mathfrak h}
(\M\vert_{\mathfrak h})\subseteq\mathcal V_{\mathfrak h}(\mathsf{N}).$ Use that result for suitably
chosen restricted subalgebras $\mathfrak h$ and their modules $\mathsf N$  to compute the support
varieties for the restricted Cartan-type Lie algebras.} 
\item[{\rm (b)}] \emph{The {\rm complexity} $\csf_{\g}(\M)$ of  a $\uf(\g,\chi)$-module
$\M$ is the rate of growth of a minimal projective resolution of $\M$ and is also given by $\csf_{\g}(\M)
=\dimm\mathcal V_\g(\M)$. Compute $\csf_{\g}(\Sf)$ for the simple $\g$-modules $\Sf$.}
\end{itemize}
\end{itemize}
\smallskip

The result discussed in part (a) played an essential role in the work of Feldvoss and Nakano
\cite{FeN}, which determined $\mathcal V_\g(\Sf)$ and $\csf_{\g}(\Sf)$ for the simple modules
$\Sf$ of the Witt algebra $\witt$. For other Cartan-type Lie algebras, some results on their
support varieties can be found in \cite{YSh}.

For a restricted Lie algebra $\g$, the complexity of  a $\uf(\g,\chi)$-module $\M$ satisfies $\csf_{\g}(\M)
\leq \dimm\mathsf{Ext}_{\uf(\g,\chi)}^{2n}(\M,\M)$ for every $n\geq 1$, and so if $\mathsf{Ext}_{\uf
(\g,\chi)}^{2n}(\M,\M)=0$ for some $n\geq 1$, the module $\M$ is projective (cf.~\cite[Lem.~2.1
and Prop.~2.2\,(1)]{F1}).

As shown by Zassenhaus \cite[Thm.~7]{Za}, the indecomposable modules of a restricted Lie algebra
$\g$ have arbitrarily large dimension. However, the dimensions of the indecomposable modules for a
particular reduced enveloping algebra $\uf(\g,\chi)$ may be bounded (see \cite{FeSt} and \cite{Pr3}).
\smallskip

\begin{itemize}
\item[]\textbf{Problem 10.} \emph{Let $\g$ be a restricted simple Lie algebra of Cartan
type, and let $\chi\in\g^*$.}
\begin{itemize}
\item[{\rm (a)}] {\it Determine the maximal dimension of a simple $\uf(\g,{\chi})$-module and
                          a lower bound on the power of $p$ dividing all their dimensions.}
\item[{\rm (b)}] {\it Describe the projective indecomposable $\uf(\g,{\chi})$-modules. In particular,
                          determine the minimal dimension of such modules and the highest power of $p$ dividing
                          all their dimensions.}
\item[{\rm (c)}] {\it Determine when the indecomposable $\uf(\g,\chi)$-modules are of \hfil  \newline bounded dimension}.
\end{itemize}
\end{itemize}
\smallskip

For the classical Lie algebras $\g$, Premet \cite{Pr1} settled a long-standing conjecture of Kac and
Weisfeiler  \cite[Rem.~at the end of Sec.~3]{KaW1} giving a lower bound on the $p$-power dividing
the dimensions of the simple $\uf(\g,\chi)$-modules when $p$ is a good prime for $\g$.
Later in \cite{Pr2}, Premet gave a simplified proof of this result using support varieties  (see also \cite{BeMRu} 
for another proof when $p$ is greater than the Coxeter number and \cite{PrSk} for more general results). 
An exposition of Premet's work on
the Kac-Weisfeiler conjecture can be found in \cite{J3}, which also establishes reciprocity rules and
shows that the blocks of the reduced enveloping algebras are determined by linkage classes of an
appropriate subgroup of the Weyl group of the classical Lie algebra $\g$. For the Poisson Lie algebra
$\g$, which is a one-dimensional central extension of the Hamiltonian Lie algebra $\mathsf{H}(2n,
\underline{1})$, Skryabin \cite{Sk1} has obtained results on the lower bound of the $p$-power
dividing the dimensions of the simple $\uf(\g,\chi)$-modules for certain characters $\chi\in\g^*$.

In \cite[Sec.~1.2]{KaW1}, Kac and Weisfeiler conjectured that the maximal dimension of a simple
$\g$-module is $p^{\half(\dimm\g-\mathsf{d})}$, where $\mathsf{d}:=\mathsf{min}\{\dimm
\mathcal{Z}_\g(\chi)\mid\chi\in\g^*\}$. They observed that their conjecture is true for classical
Lie algebras and for supersolvable restricted Lie algebras. It can be verified to hold for the Witt
algebra $\witt$, and as shown in \cite{Sk1}, it is also valid for the Poisson Lie algebra. More generally,
Premet and Skryabin \cite[Thm.~4.4]{PrSk} showed that this conjecture is true for an arbitrary
restricted Lie algebra  under the assumption that the centralizer $\mathcal{Z}_\g(\xi)$ is a torus
for some $\xi\in\g^*$.

Parts (a), (b), and (c) of Problem 10 for arbitrary restricted Lie algebras, but assuming $\chi=0$,
were posed in Problems 4, 5, and 6 of \cite{H2}. When $\chi=0$, Theorem \ref{T:finite} can be
applied to give (c).

Part (b) has been accomplished for $\mathfrak{sl}_2(\mathbf k)$ by Pollack \cite{Po} (see also
\cite{Fi}, which investigates projective indecomposable modules for $\uf(\mathfrak{sl}_2(\mathbf k))$,
and \cite[7.2 and 8.4--8.6]{B}, which considers them for $\uf(\mathfrak{sl}_2(\mathbf k), \chi)_r$,
with $r>1$). As mentioned earlier, (b) has also been solved for the Witt algebra $\witt$ by Feldvoss
and Nakano \cite{FeN}. Seligman \cite{Se} has given an explicit construction of the primitive idempotents
in $\uf(\mathfrak{sl}_2(\mathbf k),\chi)$ projecting onto the projective indecomposable $\uf(\mathfrak{sl}_2
(\mathbf k),\chi)$-modules.

If $\uni$ is a unipotent subalgebra of an arbitrary restricted Lie algebra $\g$, then $p^{\dimm\uni}$
divides the dimension of every projective module (see \cite[Lem.~3.3]{FeSt}). Therefore, when the
minimal dimension of a projective indecomposable $\uf(\g,\chi)$-module equals $p^{\dimm\uni}$ for
some unipotent subalgebra $\uni$ of $\g$, then $p^{\dimm\uni}$ will also be the highest power of
$p$ dividing the dimension of every projective (indecomposable) $\uf(\g,{\chi})$-module. For all the
examples considered below,  this will be the case, and so we will not mention it explicitly each time.

\begin{remark}
It would be interesting to find an example of a restricted Lie algebra $\g$ and a character $\chi$ for
which the minimal dimension of a projective indecomposable $\uf(\g,{\chi})$-module and the highest
power of $p$ dividing the dimension of every projective indecomposable $\uf(\g,{\chi})$-module do
not coincide or to prove that this cannot happen.
\end{remark}

The observations about unipotent subalgebras can be used to solve the second part of (b) in
Problem~10 for the restricted enveloping algebra of a classical  Lie algebra in the following way:
Let $\Phi^+$ denote the set of positive roots of a classical Lie algebra $\g$. Then the sum of the
root spaces $\uni^\pm:=\bigoplus_{\alpha\in\pm\Phi^+}\g_\alpha$ of $\g$ with respect to a maximal
torus $\mathfrak t$ is a unipotent  subalgebra with $\dimm\uni^\pm=\vert\Phi^+\vert$, and $\g=
\uni^-\oplus\mathfrak t\oplus\uni^+$ affords a triangular decomposition of $\g$. It is well known
that the Steinberg module is the unique simple $\uf(\g)$-module that is projective, and its dimension
is $p^{\vert\Phi^+\vert}$. As $p^{\vert\Phi^+\vert}=p^{\dimm\uni^\pm}$, this is the minimal
dimension of the projective indecomposable $\uf(\g)$-modules. Since baby Verma modules corresponding
to a regular nilpotent character are simple,  and one of these simple baby Verma modules is projective
(see \cite [Thm.~4.2 and Thm.~4.3\,(c)]{FrP4}), a similar result holds for the projective indecomposable
$\uf(\g,\chi)$-modules for the regular nilpotent characters $\chi$ of a classical Lie algebra $\g$. 
Moreover, for $\g$ classical and $\chi=0$ or $\chi$ regular nilpotent,  the minimal dimension of the
projective indecomposable $\uf(\g,\chi)$-modules equals the maximal dimension of the simple
$\uf(\g,\chi)$-modules.

For the Witt algebra $\g=\witt$ ($p\ge 5$), the minimal dimension $m$ of the projective indecomposable
$\uf(\g,\chi)$-modules depends only on the height $\mathsf{h}(\chi)$ of the character $\chi$. More
specifically, $m$ is given by the following expressions:
$$m= \begin{cases} p^{p-2} & \quad  \hbox{\rm if}\ \ \mathsf{h}(\chi)=-1,0,1,\\
p^{p-s-1} & \quad \hbox{\rm if} \ \ 2 \leq \mathsf{h}(\chi)\leq p-1, 
\end{cases}
$$
where $s:=\lfloor\half\mathsf{h}(\chi)\rfloor$ (see \cite{FeN}). For $p\ge 5$, Nakano \cite[Thm.~2.5.9]{N}
computed the dimensions of the projective indecomposable $\uf(\g)$-modules for the Jacobson-Witt algebra $\g=\mathsf{W}(n,\underline{1})$, and the minimal dimension is $p^{\frac{n(n-1)}{2}+n(p^n-n-1)}$. For
the Hamiltonian Lie algebra $\g=\mathsf{H}(2,\underline{1})$ (resp.~the contact Lie algebra $\g=\mathsf{K}
(3,\underline{1})$) with $p\ge 5$, the minimal dimension of the projective indecomposable $\uf(\g)$-modules
is $p^{p^2-6}$ (resp.~$p^{p^3-6}$) (see \cite[Thm.~4.6]{HoN2} and \cite[Lem.~6.3 and Thm.~6.5]{Ho2}).

For any restricted Lie algebra $\g$, part (c) of Problem 10 is equivalent to determining the reduced
enveloping algebras $\uf(\g,\chi)$ that  have finite representation type. Indeed, one direction of that
assertion is trivial, and the other follows from Roiter's solution \cite{Ro} of the first  Brauer-Thrall conjecture.
Thus, (c) of Problem~10 could follow from knowing which blocks of $\uf(\g, \chi)$ have finite representation
type (see Problem 8\,(a)).
\medskip

\begin{itemize}
\item[]\textbf{Problem 11.} {\it Every Cartan-type Lie algebra has a filtration $\g=\g_{(-q)}
\supset\g_{(-q+1)}\supset\dots\supset\g_{(r)}\supsetneq 0$ with $q=1$ or $2$.}
\begin{itemize}
\item[{\rm (a)}] {\it For the restricted simple Cartan-type Lie algebras, determine the relationship
between the height of a character $\chi$ {\rm (}as in \eqref{eq:hi}{\rm )} and the $p$-powers for
$\uf(\g,\chi)$ in parts {\rm (a)} and {\rm (b)} of Problem {\rm 10}.}
\item[{\rm (b)}] {\it Simple modules for solvable restricted Lie algebras have dimension a power
of $p$. Determine how the powers of $p$ in Problem {\rm 10\,(a) and (b)} are related to the
largest dimension of a simple module for a maximal solvable subalgebra of $\g$.}
\end{itemize}
\end{itemize}
\medskip

For a restricted Lie algebra $\g$, $\mathcal{N}_p(\g)=\{x\in\g\mid x^{[p]}=0\}$ is a conical
subvariety of the nilpotent variety $\mathcal N=\{x\in\g\mid x^{[p]^e}=0\ \hbox{\rm for} \
e\gg 0\}$. It was shown in \cite{Pr2}  for a connected reductive group $\mathsf G$ with Lie
algebra $\g=\mathsf{Lie(G)}$ such that $p$ is good for $\g$ and $\mathsf{G}^{(1)}$ is simply
connected that $\mathcal V_\g(\M)\subseteq\mathcal{N}_p(\g)\cap\mathcal Z_\g(\chi)$ for
any $\uf(\g,\chi)$-module $\M$. If $\mathsf{X}$ is the direct sum of all the simple $\uf(\g,
\chi)$-modules (necessarily there are only finitely many), then $\mathcal V_\g(\mathsf X)=
\mathcal{N}_p(\g)\cap\mathcal{Z}_\g(\chi)$ holds by \cite[Thm.~2.4]{Pr3}.

Humphreys \cite[Sec.~10]{H3} has conjectured that for an arbitrary character $\chi$ of the
Lie algebra $\g=\mathsf{Lie(G)}$ there are at most $p^{\rank \g}$ nonisomorphic simple
$\uf(\g,\chi)$-modules, where $\rank\,\g$ is the dimension of a Cartan subalgebra of $\g$.
Since for more general restricted Lie algebras the dimension of a Cartan subalgebra is not an 
invariant, we make the following conjecture using instead the maximal toral rank defined in
\eqref{eq:mt}.

\begin{conj}\label{C:iso} Let $\g$ be an arbitrary restricted Lie algebra, and let $\chi$ be a
character of $\g$.
\begin{itemize}
\item[{\rm(i)}] There are at most $p^{\mathsf{mt}(\g)}$ nonisomorphic simple $\uf(\g,\chi)$-modules.
\item[{\rm (ii)}] The number of isomorphism classes of simple $\uf(\g,\chi)$-modules is at most the number
of isomorphism classes of simple $\uf(\g)$-modules.
\end{itemize}
\end{conj}
\smallskip

\begin{itemize} 
\item[]\textbf{Problem 12.}
\begin{itemize}
\item[{\rm (a)}] \emph{Determine what relationship, if any, exists between the height $\mathsf{h}(\chi)$
of a character $\chi$ and the number of isomorphism classes of simple $\uf(\g,\chi)$-modules for $\g$ a
restricted Lie algebra of Cartan type.}
\item[{\rm (b)}] \emph{Investigate the conical subvarieties of $\mathcal{N}_p(\g)\cap\mathcal{Z}_\g(\chi)$
for the restricted simple Lie algebras of Cartan type and their characters $\chi$.}
\item[{\rm (c)}] \emph{Prove Conjecture \ref{C:iso}.}
\end{itemize} 
\end{itemize}  
\smallskip

It is well known that for any classical Lie algebra $\g_0$ there are $p^{\mathsf{mt}(\g_0)}$
isomorphism classes of simple $\uf(\g_0)$-modules (see \cite{Cu}). For a restricted simple Lie algebra
$\g$ of Cartan type with homogeneous component $\g_0$ of degree 0, there is a one-to-one
correspondence between the isomorphism classes of simple $\uf(\g)$-modules and the isomorphism
classes of simple $\uf(\g_0)$-modules (see \cite[Thm.~4]{KaW1}). As $\mathsf{mt}(\g)=\mathsf{mt}
(\g_0)$, this implies that there are $p^{\mathsf{mt}(\g)}$ isomorphism classes of simple $\uf(\g)$-modules.
Thus, parts (i) and (ii) of Conjecture \ref{C:iso} are equivalent for any restricted simple Lie algebra
$\g$ of classical or Cartan type and any character $\chi$ of $\g$.

\begin{section}{Evidence for Conjecture \ref{C:iso} and a related conjecture}\end{section} 

We now present some examples that substantiate the validity of Conjecture \ref{C:iso}. 

\smallskip

\noindent $\bullet$ \ When $\g=\mathfrak{sl}_2(\mathbf k)$ and $p\geq 3$, then every nonzero character
$\chi$ of $\g$ is either regular semisimple or regular nilpotent, and it follows from \cite[Sec.~2]{FrP4}
(see also \cite[Sec.~5.2]{StF}) that the number $\vert\mathrm{Irr}(\g,\chi)\vert$ of isomorphism
classes of simple $\uf(\g,\chi)$-modules is
$$
\vert\mathrm{Irr}(\g,\chi)\vert=
\begin{cases}
\qquad p & \quad\hbox{\rm if} \ \
\chi\mbox{ is regular semisimple or }\chi=0,\\
\half(p+1) & \quad \hbox{\rm if} \ \ \chi\mbox{ is regular nilpotent}.
\end{cases}
$$
\noindent As a consequence, Conjecture \ref{C:iso} holds for $\mathfrak{sl}_2(\mathbf k)$.
\medskip

\noindent $\bullet$ \ Assume that $p$ is good for an arbitrary classical Lie algebra $\g$. Then the
$\mathfrak{sl}_2(\mathbf k)$ example can be generalized to certain characters $\chi$ of $\g$ as
follows: \ If $\chi$ is semisimple, hence  in particular if $\chi=0$, then there are $p^{\mathsf{mt}
(\g)}$ isomorphism classes of simple $\uf(\g,\chi)$-modules (see \cite[Thm.~3.5]{FrP4}). As a result,
Conjecture~\ref{C:iso} holds in this case as well.

Conjecture \ref{C:iso} is not known for arbitrary nilpotent characters of a classical Lie algebra, but in some special instances
it does hold. For example, if $\chi$ is a regular nilpotent character, then every baby Verma module
of $\uf(\g,\chi)$ is simple (see \cite[Thm.~4.2]{FrP4}); or more generally, if $\chi$ is a nilpotent
character of standard Levi form, then every baby Verma module of $\uf(\g,\chi)$ has a unique
maximal submodule (see \cite[Prop.~2.3 and Cor.~3.5]{FrP5}, \cite[Sec.~17]{H3}, and
\cite[Lem.~D.1]{J6}). Both of these results on baby Verma modules imply that $p^{\mathsf{mt}
(\g)}$ is an upper bound for the number of isomorphism classes of simple $\uf(\g,\chi)$-modules.
Since the nilpotent characters of $\mathfrak{sl}_n(\mathbf k)$ and the subregular nilpotent
characters of $\mathfrak{so}_{2n+1}(\mathbf k)$ have standard Levi form, Conjecture \ref{C:iso}
has an affirmative answer for such characters.

Assume $\g=\mathfrak{sl}_n(\mathbf k)$, and let $\chi=\chi_s+\chi_n$ be the Jordan-Chevalley
decomposition of an arbitrary character $\chi$ of $\g$. If $\chi_s=0$ (i.e., $\chi=\chi_n$ is
nilpotent), then it follows from the previous paragraph that Conjecture \ref{C:iso} holds.
Otherwise, $\chi_s\ne 0$, and then $\mathfrak{l}:=\mathcal{Z}_\g(\chi_s)$ is a Levi subalgebra
of $\g$ containing the maximal torus of $\g$ and the root spaces $\g_{\pm\alpha}$ for which
$\chi_s([\g_\alpha,\g_{-\alpha}])=0$. Since $\chi_s$ is semisimple and nonzero, there exists
at least one simple root $\alpha_\bullet$ such that $\chi_s([\g_{\alpha_\bullet},\g_{-\alpha_\bullet}])
\ne 0$ (for more details see \cite[Lemma 3.1]{FrP4} and \cite[B.8]{J6}). It follows from the
Kac-Weisfeiler Reduction Theorem (see \cite[Thm.~2]{KaW1}, \cite[Thm.~3.2 and Thm.~8.5]{FrP4},
or \cite[Thm.~B.8]{J6}) that $\vert\mathrm{Irr}(\g,\chi)\vert=\vert\mathrm{Irr}(\mathfrak{l},\chi
\vert_{\mathfrak{l}})\vert$, and we can proceed by induction. This argument shows
that Conjecture \ref{C:iso} holds for $\mathfrak{sl}_n(\mathbf k)$, hence also for $\mathfrak{gl}_n
(\mathbf k)$, since simple modules for $\mathfrak{sl}_n(\mathbf k)$ and $\mathfrak{gl}_n(\mathbf k)$
are essentially the same.
\medskip

\noindent $\bullet$ \ The Witt algebra $\witt$ for $p\ge 5$ provides further evidence to support
the conjecture. The height $\mathsf{h}(\chi)$ of a character $\chi$ for $\g=\witt$ belongs to
$\{-1,0,\dots,p-1\}$, and $\mathsf{h}(\chi)=-1$ exactly when $\chi=0$. By (\cite{Ch}, \cite{St1},
and \cite{HuS}), we have the following:
\begin{itemize}
\item[(1)] If $\chi=0$, then every simple $\uf(\g)$-module has dimension $1,p-1$, or $p$, and there
are $p$ isomorphism classes of simple $\uf(\g)$-modules.
\item[(2)] If $\chi \neq  0$, then every simple $\uf(\g,\chi)$-module has dimension $p^{s+1}$
(resp.~$p^s$) when $\mathsf{h}(\chi) \neq p-1$ (resp.~$\mathsf{h}(\chi) = p-1$), where $s=
\lfloor \half \mathsf{h}(\chi)\rfloor$. The number $\vert\mathrm{Irr}(\g,\chi)\vert$ of isomorphism
classes of simple $\uf(\g,\chi)$-modules is given by
$$
\vert\mathrm{Irr}(\g,\chi)\vert=
\begin{cases}
1 & \quad\hbox{\rm if} \ \ \mathsf{h}(\chi)\neq 0,1,p-1,\\
p-1 & \quad\hbox{\rm if} \ \ \mathsf{h}(\chi)=0,\\
p & \quad\hbox{\rm if} \ \ \mathsf{h}(\chi)=1,\\
p-1 \ \, \hbox{\rm or} \,\ p \ \ \hbox{\rm (depending on}\,\,\chi) &\quad\hbox{\rm if} \ \ \mathsf{h}(\chi)=p-1.
\end{cases}
$$
\end{itemize}
Consequently, Conjecture \ref{C:iso} holds for $\witt$ also.
\medskip

\noindent $\bullet$ When $\g$ is a Jacobson-Witt algebra $\mathsf{W}(n,\underline{1})$ for $p\geq 5$,
and $\chi$ is a character of height $\mathsf{h}(\chi) \leq 1$, then every simple $\uf(\g,\chi)$-module
is the homomorphic image of a unique induced module $\mathsf{Z}(\LL_0,\chi)$ for some simple
$\mathfrak{gl}_n(\mathbf k)$-module $\LL_0$ (see \cite[Prop.~2.2 and Thm.~4.1]{Ho5}). Under these assumptions
on $\g$ and $\chi$, Conjecture \ref{C:iso} for the reduced enveloping algebra $\uf(\g,\chi)$ follows
from knowing that Conjecture \ref{C:iso} is true for $\mathfrak{gl}_n(\mathbf k)$ (see above and also \cite[Thm.~4.2\,(1),
Thm.~4.3\,(1), and Thm.~4.4\,(1)]{Ho5}).

According to \cite[Thm.~4.3 and Thm.~4.4\,(2)]{HoZ1}, when $\mathsf{h}(\chi)=1$ and $p\ge
5$, the number of isomorphism classes of simple $\uf(\g,\chi)$-modules for any restricted simple Lie
algebra $\g$ of Cartan type coincides with the number of isomorphism classes of simple $\uf(\g_0,
\chi\vert_{\g_0})$-modules for its component $\g_0$ of degree 0, and therefore in this case, the
validity of Conjecture \ref{C:iso} for the reduced enveloping algebra $\uf(\g,\chi)$ of any restricted
simple Lie algebra $\g$ of Cartan type is a consequence of the validity of Conjecture \ref{C:iso} for
$\uf(\g_0,\chi\vert_{\g_0})$. In particular, Conjecture \ref{C:iso} holds for the special Lie algebras
$\g=\mathsf{S}(n,\underline{1})$ and $\chi \in \g^*$ with $\mathsf{h}(\chi)=1$, since it is true for
$\mathfrak{sl}_n(\mathbf k)$. For certain characters $\chi$ with $\mathsf{h}(\chi)>1$, Conjecture \ref{C:iso}
is also valid for the Jacobson-Witt algebras $\mathsf{W}(n,\underline{1})$ (see \cite[Thm.~4]{HoZ2})
and for the special Lie algebras $\mathsf{S}(n,\underline{1})$ (see \cite[Thm.~4.1]{Z2}).
\medskip

\noindent $\bullet$ \ It is shown in \cite[Prop.~1.2 and Cor.~1.3]{FeSiWe} that for any restricted Lie
algebra $\g$ and any $\uf(\g,\chi)$-module $\M$, the projective cover $\mathsf{P}(\M)$ satisfies
\begin{equation}
\dimm \mathsf{P}(\M)\leq(\dimm\M)\cdot p^{\dimm\g-\mathsf{mt}(\g)},
\end{equation}
and when equality holds, then $\mathsf{P}(\M)\cong\mathsf{Ind}_{\mathfrak t}^{\g}(\M,\chi)$
for any torus $\mathfrak t$ of maximal dimension. Moreover, when $p>3$ and $\M$ is the trivial
$\g$-module $\mathbf k$, then $\dimm\mathsf{P}(\mathbf k)=p^{\dimm\g-\mathsf{mt}(\g)}$
if and only if $\g$ is solvable (see \cite[Thm.~6.3]{FeSiWe}). This result is used in the proof of
\cite[Thm.~3.1]{FeSiWe} to conclude for any $\chi \in \g^*$ that the number of isomorphism
classes of simple $\uf(\g,\chi)$-modules is bounded above by $p^{\mathsf{mt}(\g)}$ when $\g$
is solvable. Consequently, Conjecture \ref{C:iso}\,(i) holds for every solvable restricted Lie algebra.
\medskip

\noindent $\bullet$ \ A \emph{$p$-envelope} of a Lie algebra $\g$ is a restricted Lie algebra
$\widehat\g$ and an injective Lie homomorphism $\imath:\g\rightarrow\widehat\g$ such that
$\langle\imath(\g)\rangle_p=\widehat\g$. The $p$-envelope of smallest dimension is unique
(as a Lie algebra) up to isomorphism and is referred to as the \emph{minimal $p$-envelope.}
Recently, (i) of Conjecture \ref{C:iso} was shown to hold for the minimal $p$-envelope of every
Cartan-type Lie algebra $\W(1,\underline{m})$ in characteristic 2 (see \cite[Thm.~3.5]{LaWe},
which generalizes the result for $\W(1,\underline{2})$ in \cite{FeSiWe}).
\medskip

The minimal $p$-envelope of a non-restricted simple Lie algebra of Cartan type is simple as a
restricted Lie algebra and semisimple when regarded as a Lie algebra. It is tempting to conjecture
that the following strengthening of Conjecture \ref{C:iso}\,(i) holds for non-restricted simple
Lie algebras of Cartan type.
\medskip

\begin{conj}\label{conj2}
Let $\mathfrak{X}(n,\underline{m})$ denote the minimal $p$-envelope of a simple Lie algebra
$\mathsf{X}(n,\underline m)^{(2)}$ of Cartan type for $\mathsf{X}=\mathsf{W,S,H,K}$.
There are at most $p^{\mathsf{mt}(\mathsf{X}(n,\underline{1}))}$ nonisomorphic simple
$\uf(\mathfrak{X}(n,\underline{m}),\chi)$-modules for a character $\chi$ of non-maximal height.
\end{conj}
 
For the algebras $\W(1,\underline{m})$ ($p\ge 5$) and characters $\chi$ of non-maximal height,
Conjecture~\ref{conj2} is a consequence of \cite[Thm.~3.6\,(1) and (2)]{Sh1}. However, if
the height of $\chi$ is maximal for $\W(1,\underline{m})$, then the conjecture is known to fail
(see \cite{Mi} and \cite[Thm.~3.6\,(3)]{Sh1}), although Conjecture \ref{C:iso}\,(i) holds in this case.
\medskip

\bibliographystyle{amsalpha}
 
\end{document}